\documentclass[fleqn,a4paper]{article}
\usepackage{a4wide}
\usepackage{amsmath}
\usepackage{amssymb}
\usepackage{comment}
\usepackage{hyphenat}
\usepackage{multicol}
\usepackage{tabularx}

\usepackage{xspace}
\usepackage{hyperref}

\usepackage{caption}
\usepackage{subcaption}
\usepackage{graphicx}

\newcommand{\qed}{\hfill$\Box$}

\newcommand{\LRB}[1]{\Bigl(#1\Bigr)}

\newcommand{\SET}[1]{\left\{#1\right\}}

\newcommand{\MOD}[1]{\left|#1\right|}

\newcommand{\LE}{{\leqslant}}

\newcommand{\wLE}{{\,\leqslant}\,} 
\newcommand{\wGE}{{\,\geqslant}\,} 

\newcommand{\eqdef}{\doteq}

\newcommand{\RR}{\mathbb{R}}
\newcommand{\ZZ}{\mathbb{Z}}

\newcommand{\BNB}{B\&B\xspace}

\newcommand{\xik}{x_{ik}}
\newcommand{\xjk}{x_{jk}}
\newcommand{\yijk}{y_{ijk}}
\newcommand{\etaijk}{\eta_{ijk}}
\newcommand{\zijk}{z_{ijk}}
\newcommand{\zetaijk}{\zeta_{ijk}}

\renewcommand{\IJ}{\mathrm{IJ}}
\newcommand{\IN}{{\in}}

\newcommand{\quot}[1]{``{#1}''}
\newcommand{\tiff}{{\bf iff}~}

\newcommand{\KIAE}{{HPC4/HPC5}\xspace}

\title{Packing of Circles on Square Flat Torus as Global Optimization of Mixed Integer Nonlinear problem}
\author{Sergey A. Smirnov\footnote{Center for Distributed Computing, Institute for Information Transmission Problems of the Russian Academy of Sciences,  Bolshoy Karetny per. 19-1, Moscow 127051, Russia, \texttt{sasmir@gmail.com, vv\_voloshinov@iitp.ru}} \and Vladimir V. Voloshinov\footnotemark[\value{footnote}]}

\begin{document}
\maketitle

\abstract{
The article demonstrates rather general approach to problems of discrete geometry: treat them as global optimization problems to be solved by one of general purpose solver implementing branch-and-bound algorithm (\BNB). 
This approach may be used for various types of problems, i.e. Tammes problems, Thomson problems, search of minimal potential energy of micro-clusters, etc. Here we consider a problem of densest packing of equal circles in special geometrical object, so called square flat torus $\RR^2/\ZZ^2$ with the induced metric. It is formulated as Mixed-Integer Nonlinear Problem with linear and non-convex quadratic constraints. 

The open-source \BNB-solver SCIP, \href{http://scip.zib.de}{scip.zib.de}, and its parallel implementation ParaSCIP, \href{http://ug.zib.de}{ug.zib.de}, had been used in computing experiments to find \quot{very good} approximations of optimal arrangements. The main result is a confirmation of the conjecture on optimal packing for N=9 that was published in 2012 by O. Musin and A. Nikitenko. To do that, ParaSCIP took about 2000 CPU*hours (16 hours x 128 CPUs) of cluster HPC4/HPC5, National Research Centre "Kurchatov Institute", \href{http://ckp.nrcki.ru}{ckp.nrcki.ru}.
}

\section{Introduction}
Densest packing problems appear in many areas of discrete geometry. 
Hereinafter one of these problems is considered: a kind of Tammes problem for so called square flat torus. Flat torus in a nutshell is a factor--space $\RR^2/\ZZ^2$ with metric induced by \quot{ordinary} Euclidean metric ($\ZZ^2$ - is an integer lattice in $\RR^2$). The problem of optimal packings of congruent circles into flat torus has been studied in \cite{bib:2012arXiv1212.0649M,bib:musin2016optimal}. As to practical reasons this problem relates to the problem of \quot{super resolution of images} for aero-photography and space imagery. 

\begin{figure}[h]
\centerline{\includegraphics[scale=0.16]{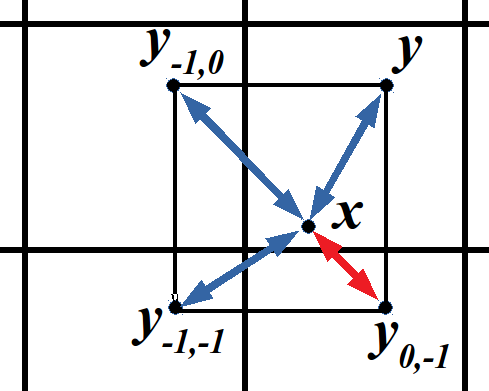}} 
\caption{Distance between $\mathbf{x}$ and $\mathbf{y}$ on Square Flat Torus}
\label{fig:ft_metric}
\end{figure}

The problem has very simple formulation: find arrangement of $N$ points in flat torus to maximize minimal distance between any pair of these points. The feature of the problem is an \quot{unusual} distance between points of a torus, which are treated as equivalence classes of factor--space $\RR^2/\ZZ^2$. Figure \ref{fig:ft_metric} illustrates definition of distance between points $\mathbf{x}$ and $\mathbf{y}$ on Square Flat Torus - the length of red segment, i.e. we have the following formula:
\begin{equation}\label{eq:ft_metric}
\begin{array}{c}
d(x,y)\eqdef \sqrt{\LRB{\min\SET{\MOD{\mathbf{x}_{1}{-}\mathbf{y}_{1}},1{-}\MOD{\mathbf{x}_{1}{-}\mathbf{y}_{1}} } }^2 + \LRB{\min\SET{\MOD{\mathbf{x}_{2}{-}\mathbf{y}_{2}},1{-}\MOD{\mathbf{x}_{2}{-}\mathbf{y}_{2}} } }^2}.
\end{array}
\end{equation}

Oleg Musin and Anton Nikitenko have studied this problem for $N{=}\SET{2,3,4,5,6,7,8,9}$ and pool results in the articles \cite{bib:2012arXiv1212.0649M,bib:musin2016optimal} Their proof of optimality is based on a
computer enumeration of irreducible contact graphs corresponding to potentially optimal arrangement of points on the flat torus. They have proved optimal packings for $N$ up to 8 and presented a conjecture of optimal arrangement for $N{=}9$. 

The approach presented hereinafter is based on global optimization by branch-and-bound solvers. The packing problem is formulated as mixed-integer nonlinear programming problem (MINLP) with binary variables and non-convex quadratic constraints. This approach is not a new one and has been widely used in studies of packing problems, e.g. see \cite{bib:castillo2008solving}. We have paid attention to flat torus packing problem in our experiments on global and discrete optimization in distributed computing environment \cite{bib:voloshinov2017implementation,bib:smirnov2017concurrent,bib:smirnov2018dd}. The main advantage of combinatorial geometry is an explicit, analytic, expression for optimal "max-min"-distance. The main advantage of global optimization -- ability to use general-purpose branch-and-bound solver including parallel implementations for high-performance clusters.

\paragraph{Paper Organisation.} This paper is organised as follows. Formulation of Flat Torus Packing problem as mixed integer nonlinear (bilinear non-convex) mathematical programming problem is done in Section \ref{sec:formulation}. The Section \ref{sec:experimentsScip} presents results of computing experiments for $N{=}\SET{4,5,6,7,8}$, including brief description of SCIP solver \cite{bib:GleixnerBastubbeEifleretal.2018} that was used for global optimization. The next Section \ref{sec:experimentsParaScip} shows results for $N=\SET{8,9}$ obtained by ParaSCIP solver \cite{bib:shinano2011parascip} that is parallel implementation of SCIP based on MPI. Some details of ParaSCIP compilation on \KIAE cluster of NRC ``Kurchatov Institute'' are provided also. Finally, we present solution found by ParaSCIP for $N{=}9$, which confirmed conjecture made in \cite{bib:2012arXiv1212.0649M, bib:musin2016optimal}. 
Conclusion Section \ref{sec:experimentsParaScip} is followed by Acknowledgements.

%
%

\section{Formulation as global optimization MINLP} \label{sec:formulation}
Let $\IJ$ be a set of unordered pairs of points' indices: $\IJ\eqdef\SET{(i,j): 1{\LE}i{<}j{\LE}N}$. The problem may be formulated as follows (maximize minimum of squared paiwise distance (\ref{eq:ft_metric})):
\begin{equation}\label{eq:tftgo}
\begin{array}{c}
D\to \max\limits_{x_{ik}}:\\
D \wLE \sum\limits_{k{=}1{:}2}\LRB{\min\SET{\MOD{x_{ik}{-}x_{jk}},1{-}\MOD{x_{ik}{-}x_{jk}} } }^2~~\LRB{(i,j)\IN\IJ},\\
0\LE\xik\LE 1 \LRB{k{=}1{:}2,~i{=}1{:}N}.
\end{array}
\end{equation}

Formulation (\ref{eq:tftgo}) has non-differentiable functions in constraints. 
Let's avoid non-smoothness at the expense of introducing auxiliary continuous and binary variables:\footnote{There are a number of literature on various (similar to each other) \quot{tricks} to perform such conversion of non-smooth problems to MILP or MINLP, but, it seems that the article \cite{bib:dantzig1960significance} was one of the first.}:
\begin{equation}\label{eq:auxvar}
\begin{array}{l}
\yijk\eqdef\min\SET{\MOD{\xik{-}\xjk},1{-}\MOD{\xik{-}\xjk} },~\LRB{k{=}1{:}2,~(i,j)\IN \IJ},\\
\zijk \eqdef {-}\MOD{\xik{-}\xjk} = \min\SET{\xjk{-}\xik, \xik{-}\xjk}~\LRB{k{=}1{:}2,~(i,j)\IN \IJ}, \\
\etaijk\IN\SET{0,1}, ~ \zetaijk\IN\SET{0,1} ~\LRB{k{=}1{:}2,~(i,j)\IN \IJ}.
\end{array}
\end{equation}

\begin{figure}[h]
\centerline{\includegraphics[scale=0.7]{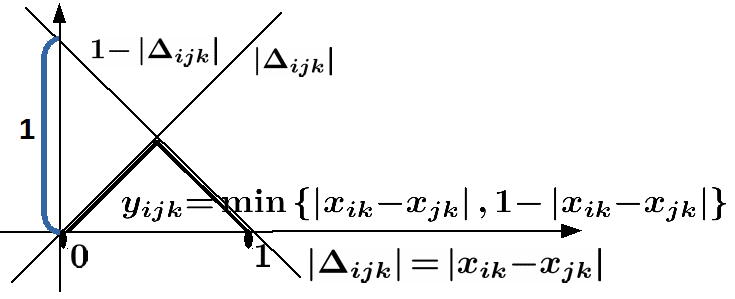}} 
\caption{Illustration for system of linear inequalities (\ref{eq:yijk})}
\label{fig:yijk}
\end{figure}

Take the first equation of (\ref{eq:auxvar}). It is equivalent to the following system of inequalities  (${k{=}1{:}2,~(i,j)\IN \IJ}$):
\begin{equation}\label{eq:yijk}
\begin{array}{l}
\yijk \wLE \MOD{\xik{-}\xjk},\\
\yijk \wLE 1{-}\MOD{\xik{-}\xjk},\\
\yijk \wGE \MOD{\xik{-}\xjk} - \etaijk\\
\yijk \wGE 1 - \MOD{\xik{-}\xjk} - 1 + \etaijk =  {-}\MOD{\xik{-}\xjk} + \etaijk.
\end{array}
\end{equation}
\noindent Equivalence means that $\yijk, \xik, \xjk$ satisfies first equation of (\ref{eq:auxvar}) \tiff there exists some binary $\etaijk$ that satisfies (\ref{eq:yijk}) with the same $\yijk, \xik, \xjk$. 
Easy proof may be done considering that 1 is the maximal difference between functions $\MOD{\Delta_{ijk}}$ and $1{-}\MOD{\Delta_{ijk}}$ on the interval $\MOD{\Delta_{ijk}}\IN[0,1]$  (inclusion $\MOD{\Delta_{ijk}}\IN[0,1]$ follows from inclusions $\xik\IN[0,1]$ and $\xjk\IN[0,1]$), see Figure~\ref{fig:yijk}.

Continue transformation of (\ref{eq:yijk}). Note that the second and the third inequalities are equivalent to the following:
\begin{equation}\label{eq:yijk23}
\begin{array}{l}
\xik{-}\xjk\wLE 1{-}\yijk, ~\xjk{-}\xik\wLE 1{-}\yijk,\\
\xik{-}\xjk\wLE \yijk{+}\etaijk, ~\xjk{-}\xik\wLE \yijk{+}\etaijk.\\
\mbox{or as two side inequalities }\\
{-}\yijk{-}\etaijk\wLE\xik{-}\xjk\wLE 1{-}\yijk,\\
{-}1{+}\yijk\wLE\xik{-}\xjk\wLE \yijk{+}\etaijk
\end{array}
\end{equation}

Let's use definition of $\zijk$ in (\ref{eq:auxvar}) to transform the first and the last lines in \ref{eq:yijk}:
\begin{equation}\label{eq:yijk14}
\begin{array}{l}
\zijk\wLE{-}\yijk,\\
\yijk\wGE\zijk{+}\etaijk.\\
\mbox{or as two side inequalities }\\
\zijk{+}\etaijk \wLE \yijk \wLE -\zijk
\end{array}
\end{equation}

\begin{figure}[h]
\centerline{\includegraphics[scale=0.7]{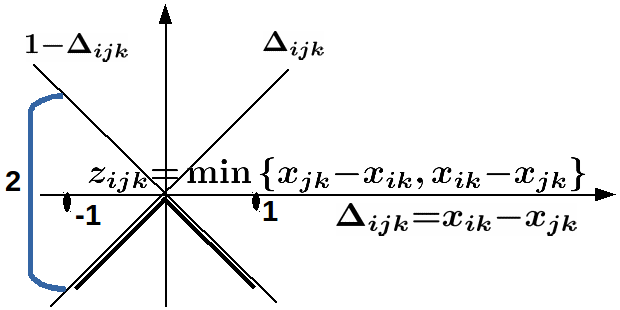}} 
\caption{Illustration for system of linear inequalities (\ref{eq:zijk})}
\label{fig:zijk}
\end{figure}

Note that definition of $\zijk$ is equivalent to the following system of linear inequalities:
\begin{equation}\label{eq:zijk}
\begin{array}{l}
\zijk\wLE\xik{-}\xjk, ~\zijk\wLE\xjk{-}\xik,\\
\zijk\wGE\xik{-}\xjk{-}2\zetaijk, \zijk\wGE\xjk{-}\xik{-}2\LRB{1{-}\zetaijk}.\\
\mbox{or as two side inequalities }\\
\zijk\wLE \xik{-}\xjk \wLE \zijk {+} 2\zetaijk,\\
{-}\zijk{-}2\LRB{1{-}\zetaijk}\wLE\xik{-}\xjk \wLE {-}\zijk.
\end{array}
\end{equation}
\noindent The proof of that equivalence is the same as mentioned after system (\ref{eq:yijk}) considering that 2 is the maximal difference between functions $\Delta$ and ${-}\Delta$ on the interval $\Delta\IN [-1,1]$, see Figure~\ref{fig:zijk}.

Finally, from definitions (\ref{eq:auxvar}) and relations (\ref{eq:yijk})--(\ref{eq:zijk}) it follows that the problem (\ref{eq:tftgo}) is equivalent to the following mixed-integer non-linear (and non-convex) problem with quadratic constraints ($k{=}1{:}2$, $i{=}1{:}N$, $(i,j)\IN \IJ$):
\begin{equation}\label{eq:minlp}
\begin{array}{l}
D\to \max ({\mbox{with variables}~\xik,\yijk,\zijk,\etaijk,\zetaijk}), \mbox{s.t.}:\\
D \wLE \sum\limits_{k{=}1{:}2}\yijk^2,\\
{-}\yijk{-}\etaijk\wLE\xik{-}\xjk\wLE 1{-}\yijk,             \\
{-}1{+}\yijk\wLE\xik{-}\xjk\wLE \yijk{+}\etaijk,             \\
\zijk{+}\etaijk \wLE \yijk \wLE -\zijk,                      \\
\zijk\wLE \xik{-}\xjk \wLE \zijk {+} 2\zetaijk,              \\
{-}\zijk{-}2\LRB{1{-}\zetaijk}\wLE\xik{-}\xjk \wLE {-}\zijk, \\
0\LE\xik\LE 1, \yijk\IN\RR, \zijk\IN\RR, \etaijk\IN\SET{0,1}, \zetaijk\IN\SET{0,1}.
\end{array}
\end{equation}
\noindent The problem has $2N^2$ continuous variables $\xik, \yijk, \zijk$ and $2N(N{-}1)$ binary variables $\etaijk, \zetaijk$. Then, except the last line, the problem has $5N(N{-}1)$ linear constraints with continuous variables, just as much linear mixed-integer constraints with continuous and binary variables and $\frac{N(N-1)}{2}$ quadratic non-convex constraints with continuous variables.

In addition to above constraints, some auxiliary constraints may be added to reduce the number of redundant solutions that might be obtained by translation by axis \quot{OX}, \quot{OY}, renumbering of points, mirror image etc.:
\begin{equation}
\begin{array}{l}
x_{11}{=}0.5, x_{12}{=}0~(\mbox{the first point is fixed to $(0.5,0)$}),\\
x_{(i{+}1)2}\wLE x_{i2}~(1\wLE i\wLE N{-}1)~(\mbox{ascending 2nd coordinates}),\\
x_{21}\wLE x_{11}.
\end{array}
\label{eq:auxcons}
\end{equation}
Simplest additional constraint, in the last line of (\ref{eq:auxcons}), reduces computing time almost twice (due to twice less volume of domain in multi-dimension space of continuous variables).

\section{Computing experiments with SCIP, {N$\mathbf{\LE}$8}} \label{sec:experimentsScip}
Computing experiments were performed with SCIP ({\bf S}olving {\bf C}onstrained {\bf I}nteger {\bf P}ro\-gram\-ming), \cite{bib:GleixnerBastubbeEifleretal.2018}. 
This is a fairly popular an open-source solver, which can be used freely for research and educational purposes.
On the home page \href{http://scip.zib.de/}{scip.zib.de} one can read: \quot{SCIP is a framework for Constraint Integer Programming oriented towards the needs of mathematical programming experts \dots ~as a pure MIP and MINLP solver or as a framework for branch-cut-and-price.
SCIP is implemented as C callable library and provides C++ wrapper classes for user plugins. It can also be used as a standalone program to solve mixed integer programs given in various formats such as MPS, LP, flatzinc, CNF, OPB, WBO, PIP, etc.} 

In our studies on optimization modelling we prefer another, so called NL-format (representing an instance of mathematical programming problem) from AMPL (A Mathematical Programming Language) \cite{bib:ampl-book}, which has almost 35 years long history (since 1985). SCIP supports NL-format by special \texttt{scipampl} build. Originally, AMPL required usage of special commercially licensed translator: to create NL-file that might be passed to any AMPL-compatible solvers, including free ones; to parse solution SOL-files returned from solvers. But in 2005 AMPL developers disclosed internal formats of NL-files \cite{bib:gay2005writing}. Thanks to that, Pyomo ({\bf PY}thon {\bf O}ptimization {\bf M}odeling {\bf O}bjects, \href{http://pyomo.org}{pyomo.org}, \cite{bib:hart2017pyomo}, open-source and free optimization modeling tool) now supports creation of NL-files. 
Thus, very popular in scientific researches Python programming language may be used to generate optimization problems, which may be solved by a proper AMPL-compatible solver.

Important feature of SCIP is its capability to solve optimization problems having polynomials in constraints (other non-linearities are admitted also). In our experiments, the Thomson problem, which has been formulated as NLP with polynomials of the 4th degree in equality-constraints, was solved by SCIP~\cite{bib:smirnov2018dd}. Details of implementations of branch-and-bound algorithm in SCIP for the case of bilinear and non-convex polynomial constraints may be found in the article \cite{bib:vigerske2017scip}. For brevity we give the following citation:
\quot{SCIP uses convex envelopes for well-known univariate functions, linearization cuts for
convex constraints, and the classical McCormick relaxation of bilinear terms. All of these are
dynamically separated, for pure NLPs also at a solution of the NLP relaxation...}.

Returning to Flat Torus Packing Problems, Pyomo is used to create NL-file from MINLP presented in (\ref{eq:minlp}) and (\ref{eq:auxcons}). Then NL-file is processed by \texttt{scipampl} application, which returns SOL-file with solution and some LOG-file with auxiliary information about solving process. Finally, SOL-file is processed by Python code via Pyomo package features to analyse solution obtained (including creation of all illustrations presented below).

Solution times for cases $N{=}4,5,6,7,8$ are presented in Table \ref{tbl:soltimes4_8}. Problems with $N{=}4,5,6,7$ has been solved on desktop 
[CPU=Intel~Core~i7-6700~@~3.40GHz, Mem=32Gb]. For $N{=}8$ the problem has been solved on standalone server [CPU=2$\times$Xeon5620~@~2.4Ghz, Mem=32Gb]. SCIP had been run with default settings, except: relative gap had been set to $1.e-6$ and memory limit - to 28Gb (actually the worst case $N{=}8$ occupied about 8Gb).

\begin{table}[ht]
\begin{center}
\begin{tabularx}{10cm}{|l|X|X|X|X|c|}
\hline
N                 & 4 & 5  & 6   & 7    & 8     \\ \hline
Solving time, sec & 3 & 30 & 118 & 2552 & 27240 (454 min) \\ \hline
\end{tabularx}
\end{center}
\caption{Solving times for $N{=}\SET{4,5,6,7,8}$, one SCIP process}
\label{tbl:soltimes4_8}
\end{table}


%
%
%

\paragraph{Results for N=7.} The case $N{=}7$ deserves more attention as, see \cite{bib:2012arXiv1212.0649M,bib:musin2016optimal}, there are three different (up to isometric transformation) optimal arrangements, see \cite{bib:musin2016optimal}, Fig.3b--Fig.3d. \hspace{-0.2em}A few efforts have been done to find all these configurations in results of SCIP solving.

\begin{figure}[!h]
\centerline{\includegraphics[scale=0.45]{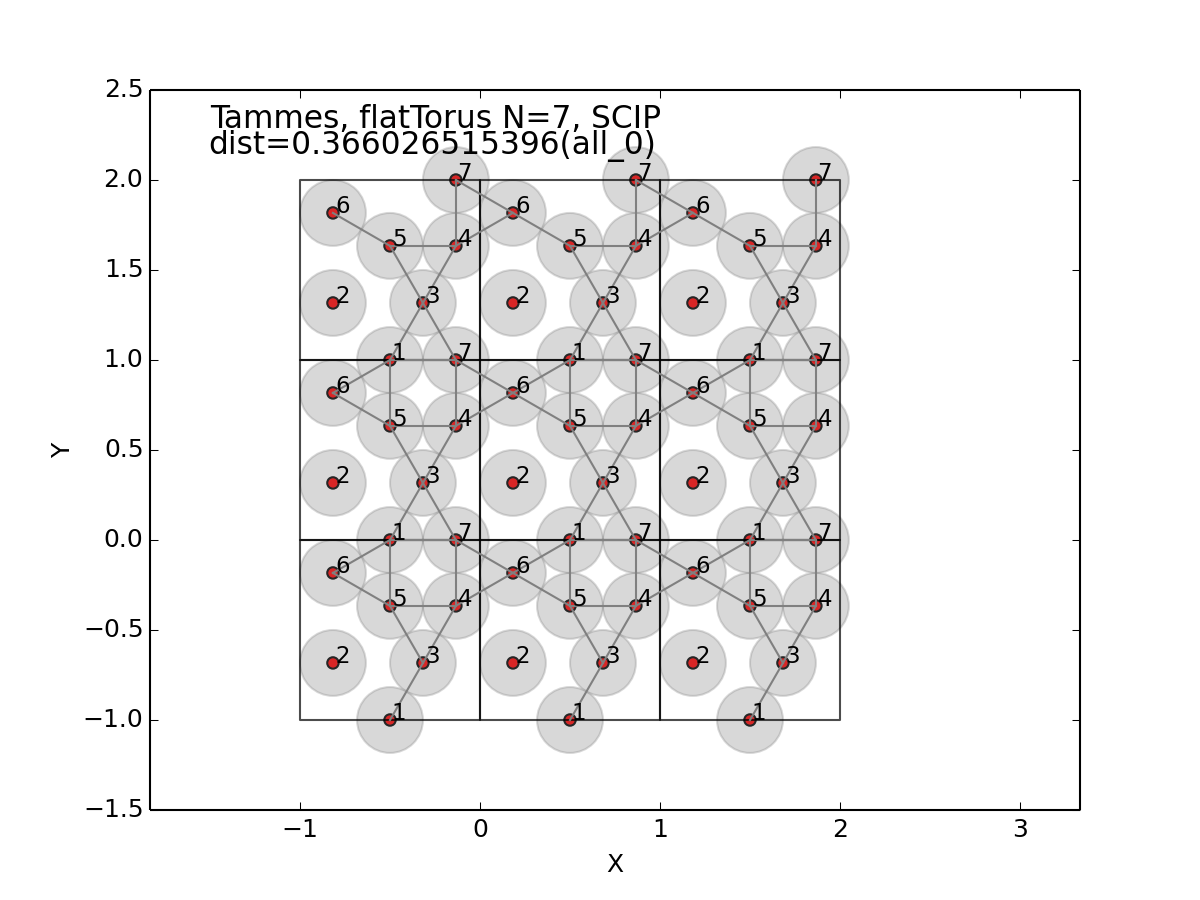}} 
\caption{$N{=}7$, the 1st configuration (see \cite{bib:musin2016optimal}, Fig.3b, $d^{*}{=}\frac{1}{1{+}\sqrt{3}}$)}
\label{fig:opt7b}
\end{figure}

\begin{figure}[!h]
\centerline{\includegraphics[scale=0.45]{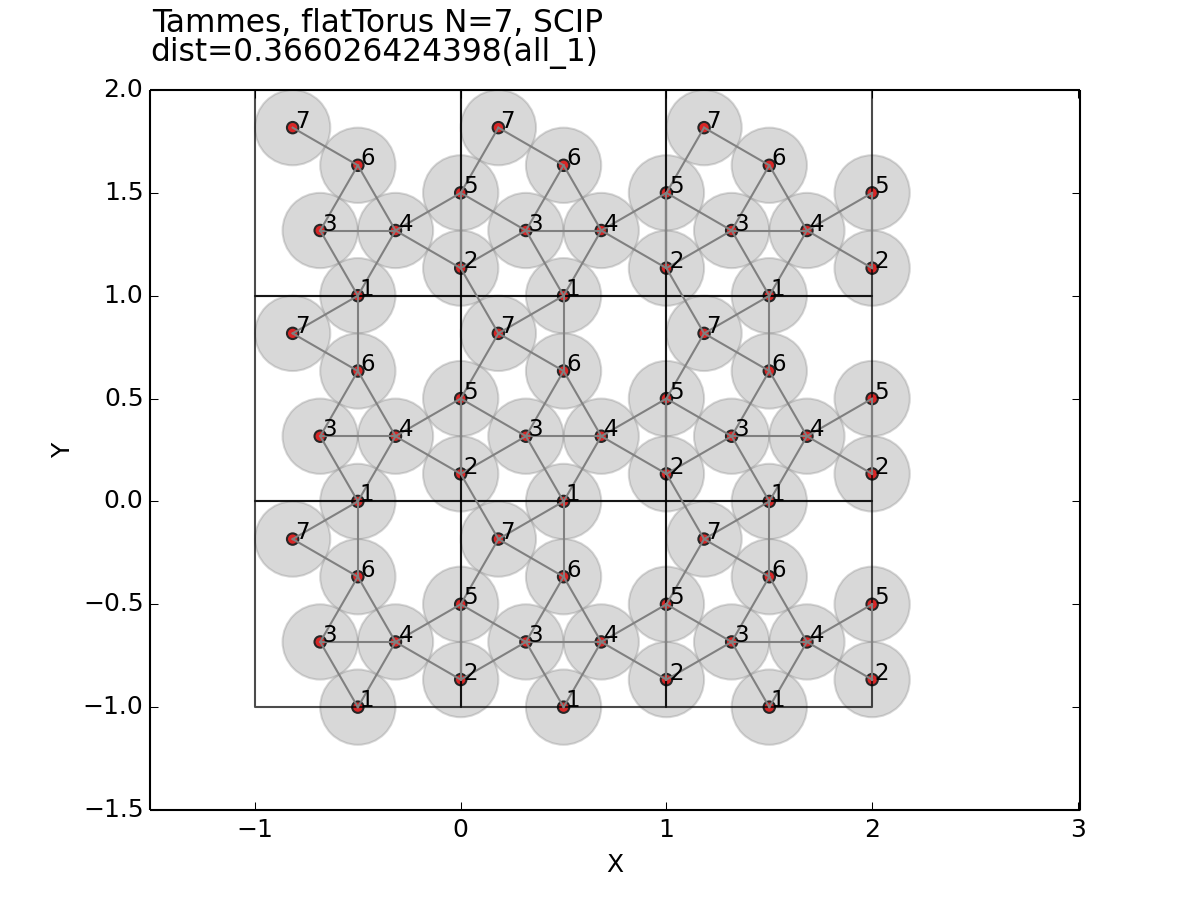}} 
\caption{$N{=}7$, the 2nd configuration (see \cite{bib:musin2016optimal}, Fig.3c, rotate $90^o \circlearrowright$ , $d^{*}{=}\frac{1}{1{+}\sqrt{3}}$)}
\label{fig:opt7c}
\end{figure}

\begin{figure}[!h]
\centerline{\includegraphics[scale=0.45]{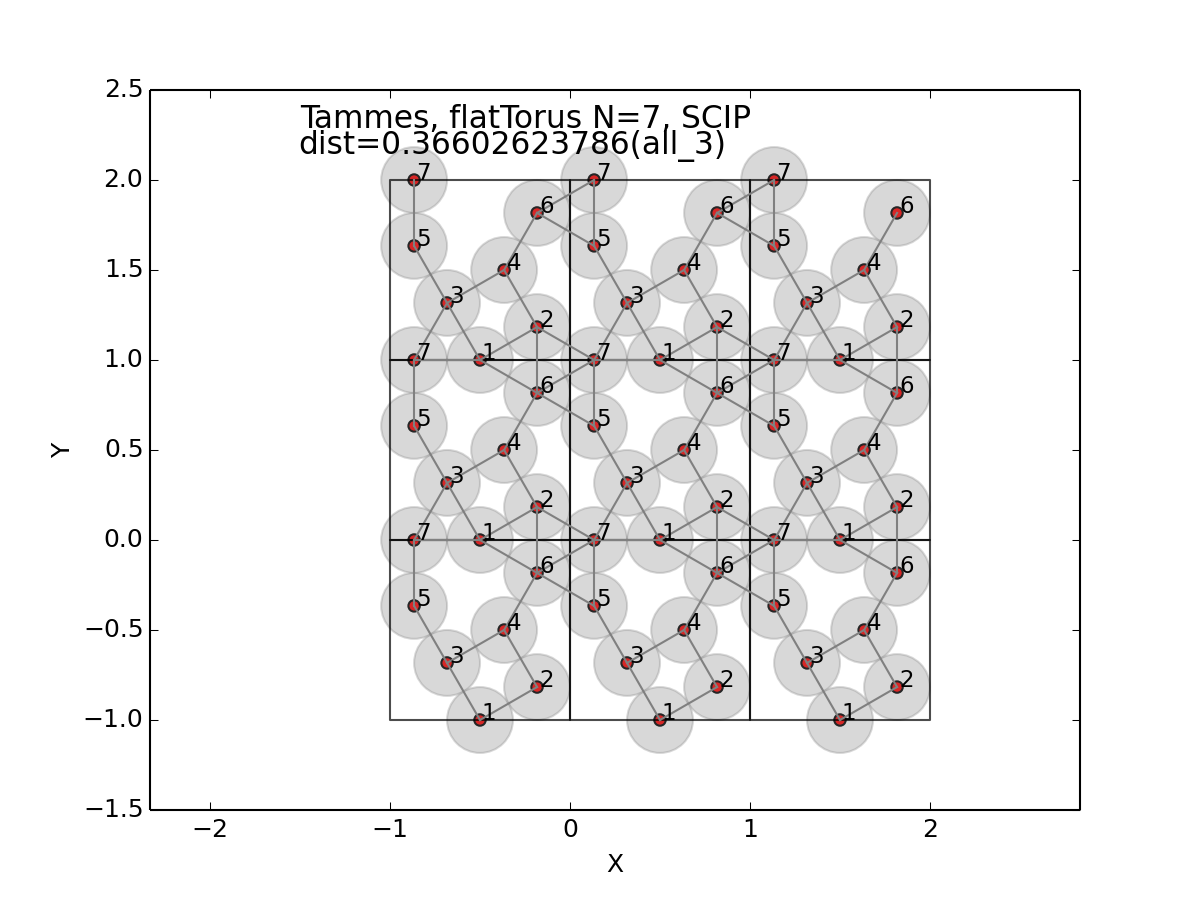}} 
\caption{$N{=}7$, the 3d configuration (see \cite{bib:musin2016optimal}, Fig.3d, flip $\updownarrow$ and rotate $90^o \circlearrowleft$ , $d^{*}{=}\frac{1}{1{+}\sqrt{3}}$)}
\label{fig:opt7d}
\end{figure}

By default SCIP solver stores optimal solutions in a list that is available by proper commands of SCIP-console. So, after successful completion of solving user can compare a number of solutions. 
Results are presented in the figures \ref{fig:opt7b}--\ref{fig:opt7d}. They have been selected manually (\quot{draw-and-compare}) from 11 optimal solutions found by standalone SCIP process. Because branch-and-bound algorithm inherently differs from \quot{exact} combinatorial method (enumeration of irreducible contact graphs) that has been used in \cite{bib:2012arXiv1212.0649M,bib:musin2016optimal}, SCIP founds redundant solutions (auxiliary constraints (\ref{eq:auxcons}) are not enough to avoid them). 

Every configurations coincides with one of those found in \cite{bib:musin2016optimal} after some isometric transformation (see captions of the figures \ref{fig:opt7b}--\ref{fig:opt7d}). Pay attention to a free position of the point number 2 on the Fig.\ref{fig:opt7b}, its circle can be freely moved within area surrounded by other grey circles.

\begin{figure}[h]\vspace*{4pt}
\centerline{\includegraphics[scale=0.5]{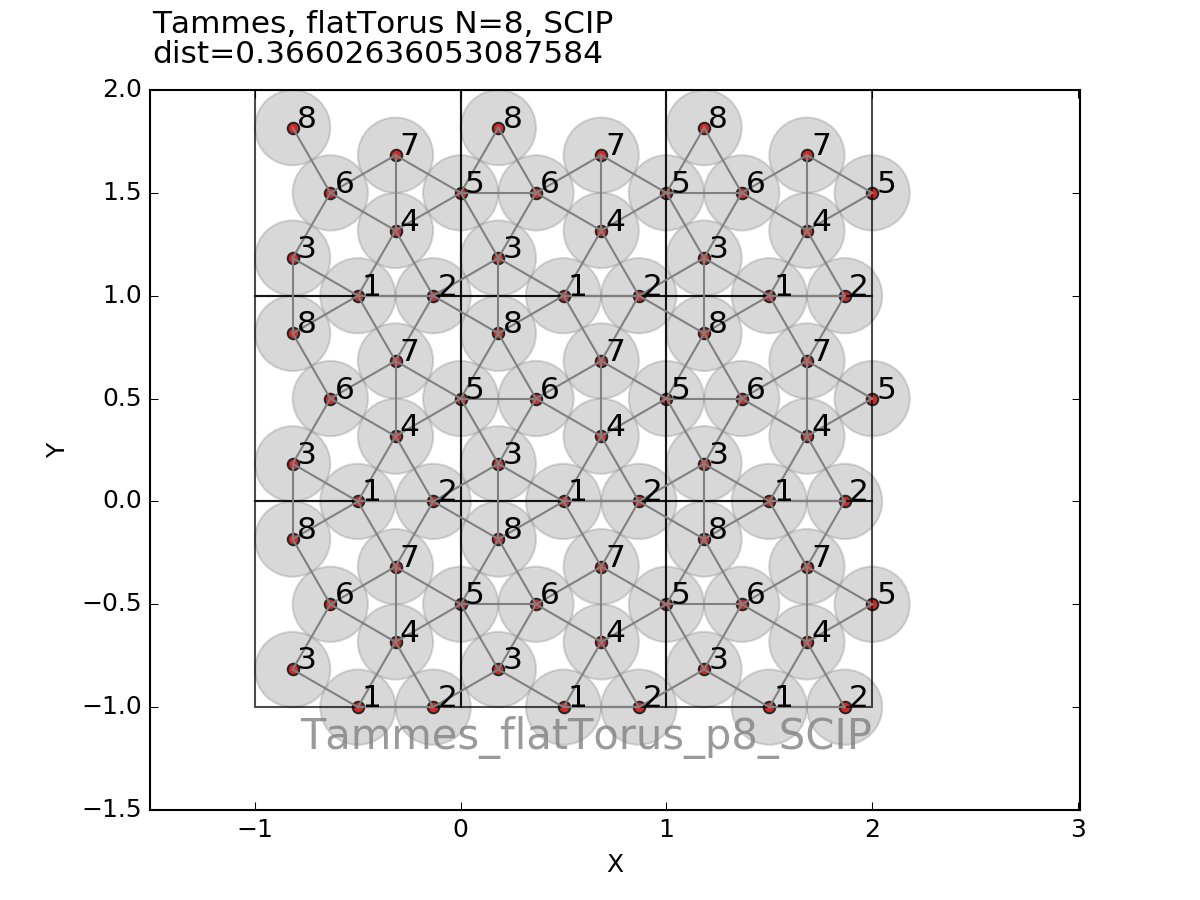}} 
\caption{Optimal configuration for $N{=}8$ found by SCIP (see \cite{bib:musin2016optimal}, Fig.3e, $d^{*}{=}\frac{1}{1{+}\sqrt{3}}$)}
\label{fig:opt8}
\end{figure}

\section{Computing experiments with ParaSCIP, {N{=}8,9}.}
\label{sec:experimentsParaScip}
The Table \ref{tbl:soltimes4_8} demonstrates dramatic growth of solving time (and complexity of Flat Torus Packing Problem) on increasing value of $N$. 
May be it is the reason for that the article \cite{bib:musin2016optimal} presents only conjecture about optimal configuration for $N{=}9$. Our attempts to solve problem for $N{=}9$ by \quot{single--threaded} SCIP have failed. 
Computing efficiency of branch-and-bound algorithm implemented in SCIP may be substantially increased by its parallel implementation as ParaSCIP solver does.

\paragraph{ParaSCIP -- parallel implementation of B\&B.}
ParaSCIP, \cite{bib:shinano2011parascip}, is a distributed memory massively parallel MIP and MINLP solver based on Ubiquity Generator (UG) framework, \href{http://ug.zib.de}{ug.zib.de}. This framework is aimed at making MIP and other discrete problem solvers parallel. In the framework the solver and the communication mechanism are abstracted making it easier to parallelize different solvers. Two communication mechanisms are available in the library: distributed memory MPI and shared memory POSIX-threads. SCIP, CPLEX and Xpress solvers are supported with MPI and only SCIP with both mechanisms. Only SCIP-based specializations of the framework (ParaSCIP for MPI-based and FiberSCIP for POSIX-threads based implementations) are publicly available. ParaSCIP can utilize quite big computing resources which allows it to solve problems not solvable even with commercial solvers on a single host: ParaSCIP ran on 80000 cores in parallel while solving open MIP instances from MIPLIB2010 \cite{bib:shinano2016solving}. 


\begin{sloppypar}
Compiling SCIP and ParaSCIP on \KIAE (cluster of National Research Center "Kurchatov Institute", \href{https://www.top500.org/site/50615}{www.top500.org/site/50615}) is not easy because \KIAE has CentOS 6 with GCC 4.4 tool-chain installed which does not support C++11 extensions required for SCIP build. We used another machine with a similar CentOS version but with devtoolset-7 (GCC 7.3) installed to build solvers and then copied them to the cluster. Machines had Intel CPUs of different families (Haswell on cluster vs westmere on our machine) so default optimizations during compilation could be not optimal for the cluster. Finally Ipopt, SCIP and ParaSCIP were compiled on our machine and then copied to the cluster where they were running well.
\end{sloppypar}

After consultation with ParaSCIP developers we ran \verb|parascip| with non-default settings. We used the following ParaSCIP settings:\\
\begin{verbatim}
Quiet = FALSE                
LogSolvingStatusFilePath = "./logs/"
LogNodesTransferFilePath = "./logs/"
CheckpointFilePath = "./logs/"
SolutionFilePath = "./logs/"
NotificationInterval = 10
LogSolvingStatus = TRUE
Checkpoint = TRUE
CheckpointInterval = 1800
\end{verbatim}

Also \quot{depth-first search} option were set for search tree traversing to reduce memory consumption in solvers:
\begin{verbatim}
nodeselection/dfs/stdpriority = 300000
\end{verbatim}

The following command line was used to run ParaSCIP:
\begin{verbatim}
mpirun parascip parascip.set tammesTorus_d2_p9.cip -q -s dfs.set
\end{verbatim}
Here \verb|parascip.set| is the ParaSCIP settings file listed above; \verb|tammesTorus_d2_p9.cip| is the problem definition converted to CIP-format from NL-format; \verb|dfs.set| is the settings file for depth-first search specified above.

\paragraph{ParaSCIP vs SCIP for N=8.} This case have been used to compare performance of ParaSCIP and SCIP on \KIAE. Solving times are presented in the Table \ref{tbl:soltimes8_KIAE}. 
One can pay attention that solving with a single SCIP process on the cluster took much more time than that on a standalone server (see last column of the Table \ref{tbl:soltimes4_8}). The reason is that the problem instance passed to the cluster did not have the last of auxiliary constraints (\ref{eq:auxcons}).
\begin{table}[!ht]
\begin{center}
\begin{tabular}{|l|l|c|} 
\hline
\parbox[t]{6cm}{\KIAE, \\ CPU Xeon E5-2680v3 12C @ 2.5GHz \vspace{0.2em}} & cores         & \parbox[t]{2cm}{solving time, \\min } \\ \hline %
SCIP      & 1             & 780               \\ \hline
ParaSCIP  & 8 (7 solvers) & 126               \\ \hline
\end{tabular}
\end{center}
\caption{Solving times for $N{=}8$, for SCIP and ParaSCIP (on \KIAE)}
\label{tbl:soltimes8_KIAE}
\end{table}

To evaluate efficiency of parallelization one should remember one feature of ParaSCIP: one of ParaSCIP processes (dedicated for MPI application), plays role of Load Coordinator and, actually, does not work with branch-and-bound algorithm's search tree. So we have:\\
efficiency (CPU): 780/126/8 = 0.77\\
efficiency (solvers): 780/126/7 = 0.88.

%

\paragraph{Results of ParaSCIP for N=9.} 
ParaSCIP with 128 processes on 8 nodes (127~solvers) of cluster \KIAE had solved the problem in 956 minutes. 

\begin{figure}[!h]\vspace*{4pt}
\centerline{\includegraphics[scale=0.5]{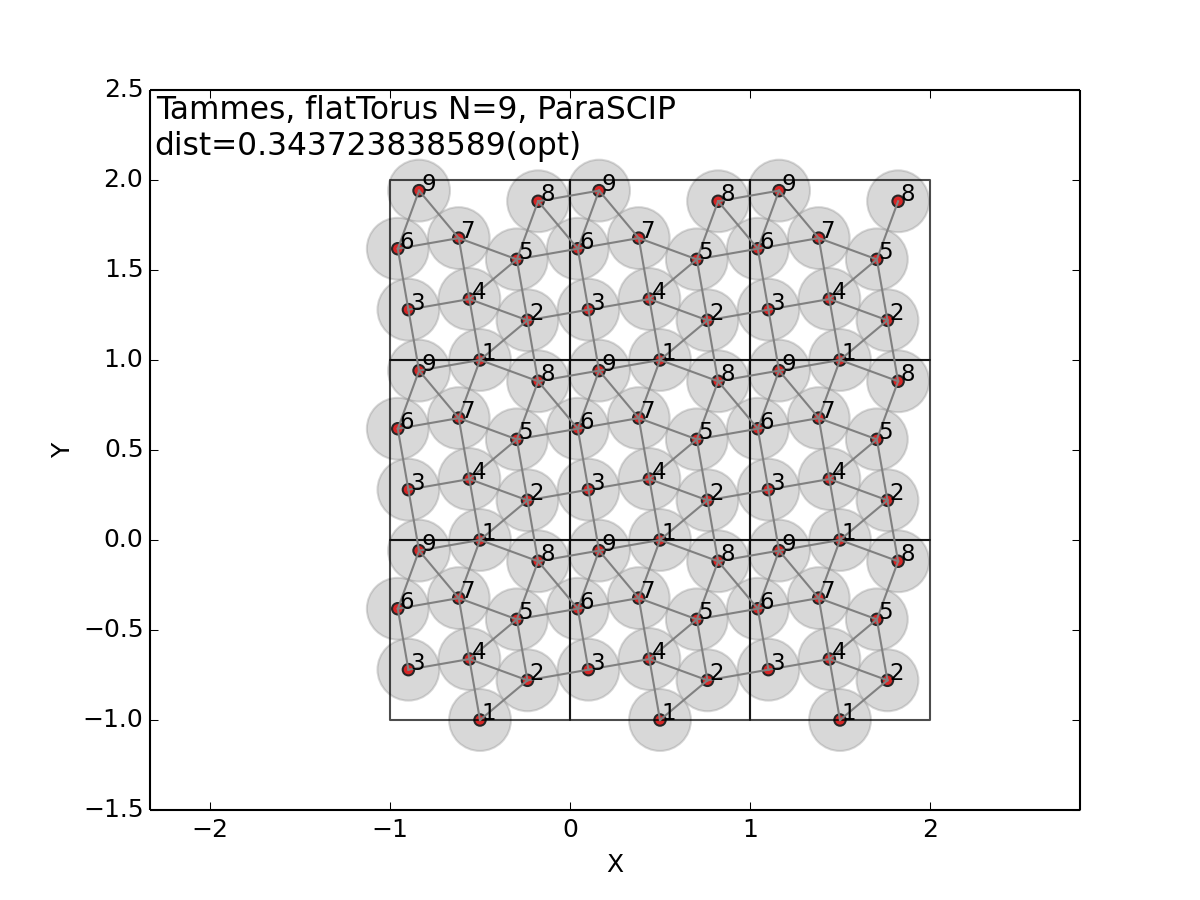}} 
\caption{Optimal configuration for $N{=}9$ found by ParaSCIP (see \cite{bib:musin2016optimal}, Fig.3f, $d^{*}{=}\frac{1}{\sqrt{5{+}2\sqrt{3}}}$)}
\label{fig:opt9}
\end{figure}

It is the most promising result presented in this work. Optimal configuration found is shown in the Fig. \ref{fig:opt9} and coincides with conjecture presented in \cite{bib:2012arXiv1212.0649M,bib:musin2016optimal}. Taking into account load balancing between working processes the following evaluation of complexity may be done: $127{\times}956/60{\approx}2043 ~\mbox{CPU}{\times}\mbox{hours}$. Total complexity, including CPU dedicated for Load Coordination is $128{\times}956/60{\approx}2059 ~\mbox{CPU}{\times}\mbox{hours}$. We believe that this complexity may be reduced almost twice if the last of auxiliary constraints (\ref{eq:auxcons}) would be added.

%

\section{Conclusion} \label{sec:conclusion}
Obtained results confirmed the relevance of global optimization approach in solving hard problems of combinatorial geometry. A few tricks in formulation Flat Torus Packing problem as mixed integer nonlinear problem and \quot{general purpose} solver SCIP with its parallel version ParaSCIP let to give \quot{computer aided} proof of optimal arrangement for 9 points. All results obtained before for $N{=}4,5,6,7,8$ have been reproduced also. 

Once more has been confirmed semi-empirical rule (see \cite{bib:smirnov2018dd}) that simple auxiliary constraints that reduces the volume of feasible domain (i.e. last inequality of (\ref{eq:auxcons})) may substantially reduce solving time for branch-and-bound algorithm using McCormik envelopes for global optimization of nonlinear problems.

General purpose solver SCIP and ParaSCIP, its parallel implementation, can be used successfully for global optimization of mixed integer nonlinear problems with bilinear constraints.

\subsubsection*{Acknowledgements}
Authors are grateful to Alexey Tarasov for recommendation to try Flat Torus Packing problem in our experiments with distributed implementations of \BNB. 
Authors thank Stefan Vigerske and Yuji Shinano for their consultations on using of SCIP and ParaSCIP solvers.\\
This work is supported by the Russian Science Foundation (project No. 16-11-10352).\\
This work has been carried out using computing resources of the
federal collective usage centre Complex for Simulation and Data
Processing for Mega-science Facilities at NRC "Kurchatov Institute", \href{http://ckp.nrcki.ru}{ckp.nrcki.ru}.

\label{sect:bib}
\bibliographystyle{unsrt}
\bibliography{tammestorus}
\end{document}